%% file: JXM_ARXIV.tex
\documentclass[12pt]{article}

\usepackage{latexsym,amssymb,amsmath}
\usepackage{amsthm, amstext}
\usepackage{dsfont}
\usepackage{array, amsfonts, mathrsfs}
\usepackage{hyperref}
\hypersetup{colorlinks=true}
\hypersetup{linkcolor=blue}
\usepackage{color}
\usepackage{adjustbox}
\usepackage[T2A]{fontenc}
\usepackage[utf8]{inputenc}
\usepackage{csquotes}
\usepackage[russian,english]{babel}
\usepackage{array}
\usepackage{hhline}
\usepackage{stackrel}
\usepackage{subcaption}

\usepackage[
    backend=biber,
    natbib=true,
    url=true,
    doi=true,
    eprint=true,
      sorting=nyt,
]{biblatex}

\newcommand{\myskip}[1]{}

\newcommand{\vpad}[1]{\adjustbox{padding=0mm 1mm 0mm 1mm}{#1}}
\newcommand{\A}{\mathfrak{A}}
\newcommand{\D}{\mathfrak{D}}
\newcommand{\F}{\mathfrak{F}}
\newcommand{\G}{\alpha}
\newcommand{\M}{\mathfrak{M}}

\newcommand{\myi}{\mathrm{i}}
\newcommand{\Nc}{}

\newtheorem{remark}{Remark}

\newtheorem*{remark*}{Remark}

\newtheorem{example}{Numerical example}

\newtheorem{conjecture}{Conjecture}

\numberwithin{equation}{section}
\numberwithin{remark}{section}
\numberwithin{example}{section}
\numberwithin{figure}{section}
\numberwithin{table}{section}
\numberwithin{conjecture}{section}

\begin{document}
\title{In Search of Approximate Polynomial Dependencies\\
Among the Derivatives of \\[2Mm]
the Alternating Zeta Function
}
\author{Yuri Matiyasevich}

\date{(Published in \emph{J.\ Exper.\ Math.}, 1(2):239–256, 2025)}

\maketitle

\begin{quote}
\textbf{Abstract:}
  It is well-known that the Riemann zeta function
does not satisfy any \emph{exact} polynomial
 differential equation.
Here we present numerical evidence for the existence of
\emph{approximate} polynomial dependencies between the values of
the alternating zeta function and its initial derivatives.

\quad A number of conjectures is stated.
\end{quote}


\section{Introduction}\label{intro}

The \emph{Riemann zeta function}
 can be defined via
a \emph{Dirichlet series}, namely,
\begin{equation}
\zeta(s)=\sum_{n=1}^{\infty}n^{-s}.
\label{zeta}
\end{equation}
The series converges for  $\Re(s)>1$ but the zeta
function can be analytically extended  to the whole
complex plane with the exception of the simple pole at $s=1$.

The zeta function is one of the most mysterious
mathematical objects. The celebrated
\emph{Riemann Hypothesis} has remained unproved for more than
 a century and a half.

The study of the zeta function might
be easier if
we found a suitable differential equation
satisfied by this function.
However, this is impossible.
D.\,Hilbert
\cite{Hilbert1900}
asserted that  the zeta function does not satisfy any algebraic
differential equation. A proof was
given later by
 V.~E.~E.~Stadigh \cite{Stadigh} and,
 in a more  general setting,
 by
D.~D.~Mor\-du\-khai-Boltovskoi
\cite{mord1} (see also \cite{ostr,mord2}).
These ``negative'' results were further
extended to broader
classes of differential equations (see, for example,
 recent publication
\cite{zbMATH07556828} and references there).

For real values of the argument,
the  zeta function
was studied already by
L.\,Euler. Likewise, he
considered the \emph{alternating  zeta function}
(known also as \emph{Dirichlet eta function})
\begin{equation}
\eta(s)=\sum_{n=1}^{\infty}(-1)^{n+1}n^{-s}.
\label{eta}
\end{equation}

The two functions, \eqref{zeta} and \eqref{eta},
are interconnected:
\begin{equation}\label{etazeta}
  \eta(s)=(1-2\times2^{-s})\zeta(s).
\end{equation}
The factor $1-2\times2^{-s}$ vanishes at $s=1$,
and thus  $\eta(s)$ is an entire function.
Thanks to this,  it is often
easier to work
with $\eta(s)$. Moreover, the series \eqref{eta} converges
in a larger half-plane $\Re(s)>0$.

The above mentioned ``negative'' results for the zeta function
can be transferred to the eta function.
In particular,  the eta
function cannot satisfy any
\emph{exact} polynomial differential equation.
In this paper we  demonstrate
\emph{approximate} polynomial relations
among the values  of the alternating
zeta function and its derivatives.
This relationship was found experimentally, and
so far  no ``explanation'' has been given for these
phenomena.

More precisely,
we define several series of pairs of
polynomials with integer
coefficients
\begin{equation}\label{VWDN}
\langle V_{N}(y_{0},\dots,y_{k},\dots,y_{N}),\
W_{N}(y_{0},\dots,y_{k},\dots,y_{N-1})\rangle,
\ N=1,2,\dots\ .
\end{equation}
For each such series   we state a conjecture that
\begin{equation}\label{VWlim1}
\frac{
  V_{N}(\eta({a}),\dots,
  \eta^{\langle n \rangle}(a),\dots,
  \eta^{\langle N-1 \rangle}(a))}
  {
  W_{N}(\eta({a}),\dots,
  \eta^{\langle n \rangle}(a),\dots,
  \eta^{\langle N-1 \rangle}(a))}
\stackrel[N\rightarrow \infty]{\longrightarrow}{} 1
\end{equation}
except for countably many values of~$a$.
To support such conjectures we
demonstrate certain numerical data.
They
show that already for relatively small values of $N$
the left-hand side in \eqref{VWlim1}
can be very close to~$1$.

The paper consists of two parts.
In Part \ref{part1}
we present raw numerical data
and discuss their peculiarities.
Based on our observations,
in Part \ref{part2} we  construct
polynomials \eqref{VWDN} with the expected
property~\eqref{VWlim1}.
In the Concluding remark we briefly discuss
possible extension of our observation and conjectures
to the Riemann zeta function.

Preliminary versions of this paper \cite{preprint}
and  \cite{RG} contain more numerical data
and conjectures.

\part{Numerical observations}\label{part1}

\section{An unreachable goal}\label{discov1sect}

Let us start by considering a very ambitious and
therefore unattainable goal:   \emph{for a particular
positive integer $N$, find
a linear polynomial
\begin{equation}\label{polP}
  P(y_0,y_1,\dots,y_{N-1})=
  y_0+c_1y_1+\dots+c_{N-1}y_{N-1}-b
\end{equation}
with numerical coefficients
\begin{equation}\label{bcc}
  b, c_1,\dots,c_{N-1}
\end{equation}
such that for all $a$ }
\begin{equation}\label{Pis0}
  P\big(\eta(a),\eta'(a),\dots,
  \eta^{\langle N-1 \rangle}(a)\big)=
  0.
\end{equation}

If such a polynomial $P$  existed,
 we could find its coefficients \eqref{bcc}
in the following way.
First, we select a set
$\A=\{a_0,\dots,a_{N-1}\}$ of $N$
complex numbers
in general position.
Then, we solve the linear system
 \begin{equation}\label{cPis0}
  \eta(a_j)+c_1(\A)\eta'(a_j)+
  \dots+c_{N-1}(\A)\eta^{\langle N-1 \rangle}(a_j)=b(\A),
  \quad j=0,\dots,N-1
\end{equation}
in unknowns
\begin{equation}\label{coeffs}
b(\A),\ c_1({\A}),\ \dots,\ c_{N-1}({\A}).
\end{equation}

Our goal would be achieved  if all numbers \eqref{coeffs}
were the same for all~$\A$.
As explained in the introduction,
this is not the case.

\begin{figure}[t]
\centering
\hfill
\subcaptionbox*{$\mathfrak{A}_1$}
[.5\textwidth]{\includegraphics[height=0.3\textheight]{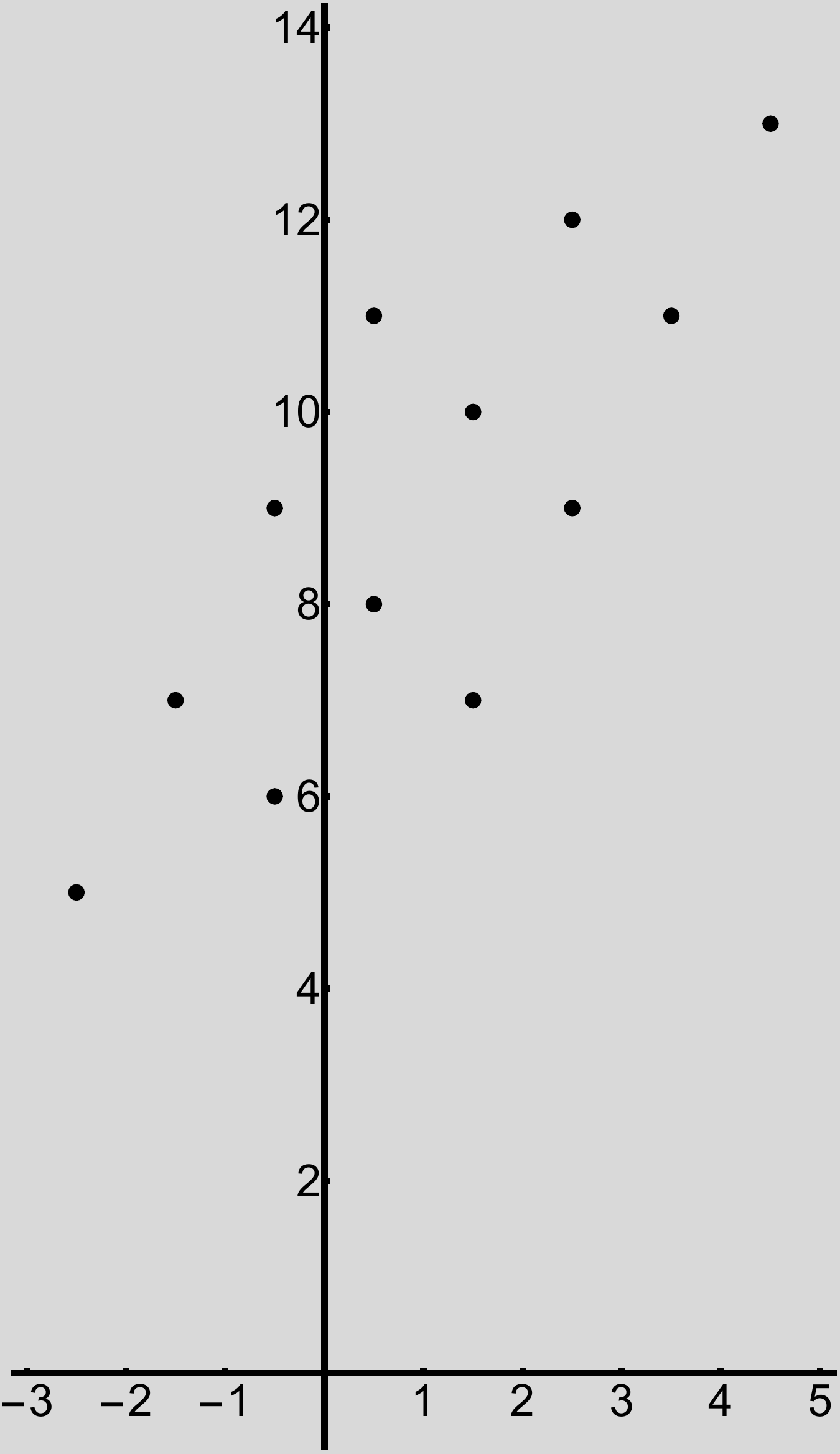}}
\hfill
\subcaptionbox*{$\mathfrak{A}_2$}
[.3\textwidth]{\includegraphics[height=0.3\textheight]{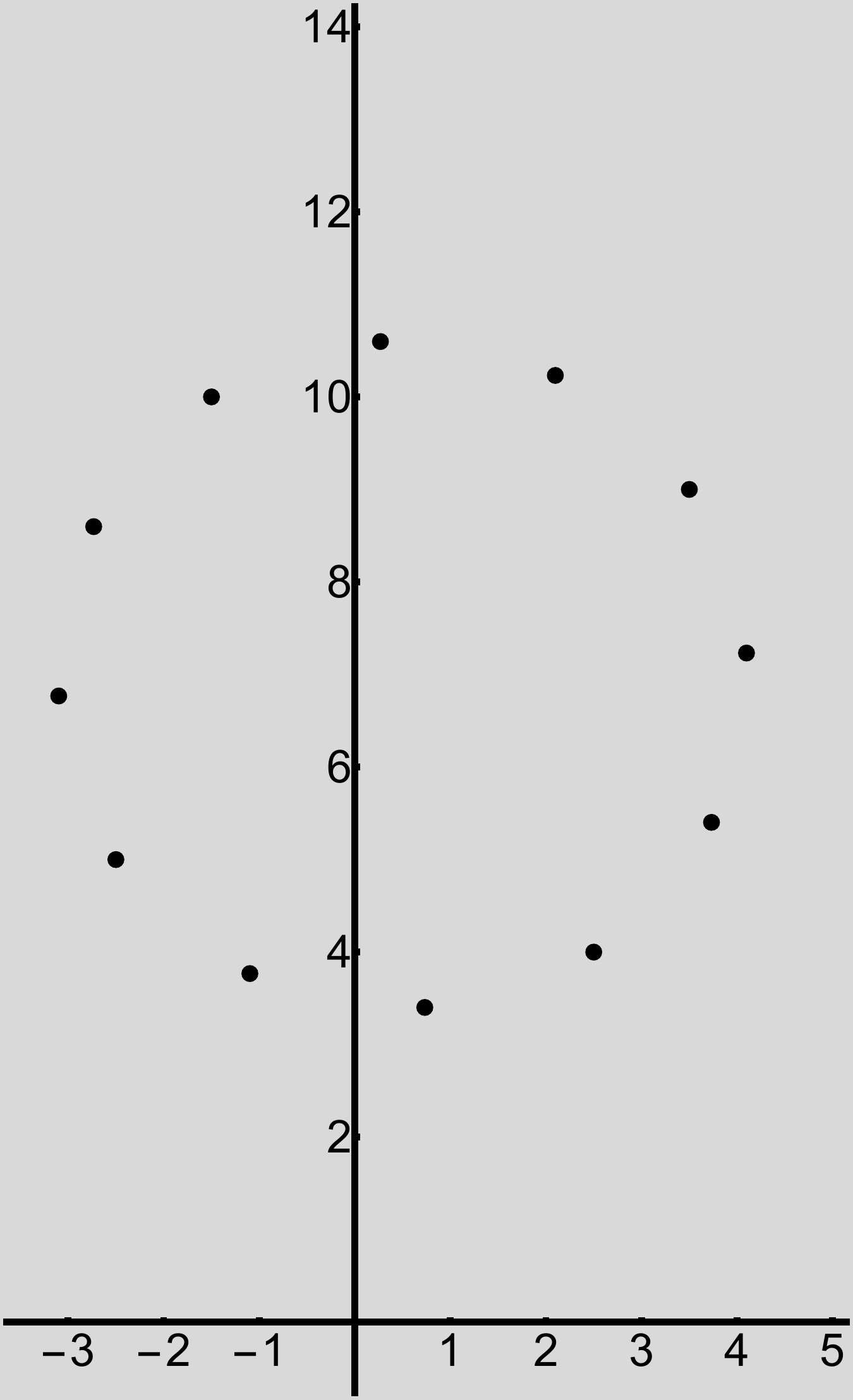}}
\hfill
\caption{Grid
$\A_1=\mathfrak{A}_\mathrm{G}(-2.5+5\mathrm{i},
2+\mathrm{i},1+2\mathrm{i},2,3)$ and  discrete circle
$\A_2=\mathfrak{A}_\mathrm{C}
(0.5+7\mathrm{i},3+2\mathrm{i},12)$
  defined respectively by
\eqref{grid} and \eqref{circle}.
}
\label{gc}
\end{figure}

Let us examine numerical data for two particular
sets depicted on Fig.~\ref{gc}. These sets are
specimens
of \emph{grids}  and \emph{discrete circles}
which in general case are defined as follows.
For complex $a$, $\delta_1$,  $\delta_2$
and non-negative integers $N_1$,~$N_2$
\begin{equation}\label{grid}
  \mathfrak{A}_\mathrm{G}(a,\delta_1,\delta_2,N_1,N_2)=\\
  \{a+k\delta_1+l\delta_2 \,\vert\,
  k=0,\dots,N_1, l=0,\dots,N_2\},
\end{equation}
and for complex $c$, $r$
and positive integer $N$
\begin{equation}\label{circle}
  \mathfrak{A}_\mathrm{C}(c,r,N)=
  \{c+\mathrm{e}^{2\pi \mathrm{i} k/N}r\, \vert\
  k=0,\dots,N-1\}.
\end{equation}

\begin{table}[t]
\centering
\input gccoeffs_2.tex
\caption{Values of the $c$-coefficients from
\eqref{cPis0}  for
 sets $\A_1$ and $\A_2$  depicted on Fig.\,\ref{gc}.}
\label{Nbc}
\end{table}

Table \ref{Nbc} presents the values of
the $c$-coefficients from the solutions of
the system \eqref{cPis0}
 for the two sets from Fig.~\ref{gc}.
We observe that while the sets are rather
different, there is certain similarity
between numbers $c_k(\mathfrak{A}_1)$
and~$c_k(\mathfrak{A}_2)$. As for the
$b$-coefficients, $b(\mathfrak{A}_1)$ and
 $b(\mathfrak{A}_2)$, they are almost equal:
 \begin{align}
   b(\mathfrak{A}_1)=
   0.99999999999999430995...
   +0.00000000000000951329...\mathrm{i},
   \label{bA1}\\
    b(\mathfrak{A}_2)=
    1.00000000000000637109...
    +0.00000000000003070321...\mathrm{i}.
    \label{ba2}
 \end{align}

\begin{table}[t]
\centering
\input multigrids.tex
\caption{Instances of accuracy of  \eqref{bis1}
for grids; $N=(N_1+1)(N_2+1)$ is the cardinality of the set
$\A=\A_\mathrm{G}(a,\delta_1,\delta_2,N_1,N_2)$.
}\label{manygrids}
\end{table}

\begin{table}[t]
\centering
\input multicircles.tex
\caption{Instances of accuracy of  \eqref{bis1}
for discrete circles.
}\label{manycircles}
\end{table}

Moreover, Tables \ref{manygrids}--\ref{manycircles}
show that  $b(\mathfrak{A})$
from the solution of
 system \eqref{cPis0}
is surprisingly close to $1$
for a selection of rather different grids and discrete circles:
\begin{equation}\label{bis1}
   b(\mathfrak{A})\approx 1.
\end{equation}
More
 relevant numerical data (for grids, discrete
circles and ``random'' sets) are presented in
\cite[Tables 36--38]{preprint}.
We return to discussion
of this phenomenon after introduction of
some notation.

\section{Matrix notation}\label{detsect}

We
rewrite
system \eqref{cPis0}  in matrix notation
and  solve it by Crammer's rule.
To this end we  introduce the following notation.

Let $D_N$ denote the determinant of the ``canonical''
$N\times N$ matrix:
\begin{equation}\label{DN}
  D_N=\begin{vmatrix}
    x_{0,0}&\dots&x_{0,N-1}\\
    \vdots&\ddots&\vdots\\
     x_{N-1,0}&\dots&x_{N-1,N-1}
  \end{vmatrix}.
\end{equation}

If $E$ is some expression and $m$ and $n$ are
particular integers such that $0\le m<N$ and
 $0\le n<N$, then
 \begin{equation}\label{repl1}
  D_N(x_{m,n}\rightarrow E)
 \end{equation}
denotes the result of replacing $x_{m,n}$ by~$E$
in~$D_N$.

Multiple replacements like
\begin{equation}\label{repl2}
  D_N(x_{1,1}\rightarrow E_1,x_{2,2}\rightarrow E_2)
\end{equation} are possible with the natural
meaning.

Several similar replacements
can be defined by a
rule with parameters $j$ and $k$.
For example,
$D_N(x_{1,k}\rightarrow E_k)$
is a short for
\begin{equation}\label{replN}
  D_N(x_{1,0}\rightarrow E_0,\ \dots,\
  x_{1,N-1}\rightarrow E_{N-1}),
\end{equation}
and $D_N(x_{j,k}\rightarrow E_{j,k})$
is a short for
\begin{multline}\label{replN2}
  D_N(x_{0,0}\rightarrow E_{0,0},\ \dots,\
  x_{0,N-1}\rightarrow E_{0,N-1},
  \dots,\\
  x_{N-1,0}\rightarrow E_{N-1,0},\ \dots,\
  x_{N-1,N-1}\rightarrow E_{N-1,N-1}).
\end{multline}

Replacement rules are applied consecutively
from left to right, for example, in
\begin{equation}
  D_N(x_{0,k}\rightarrow E_k,\
  x_{j,0}\rightarrow F_{j})
\end{equation}
$x_{0,0}$ is replaced by $E_0$.

A replacement rule can be conditional, for example,
\begin{equation}\label{C}
  C_N=D_N(j=k\Rightarrow x_{j,k}\rightarrow x_{j,k}-\lambda)
\end{equation}
 is the characteristic polynomial of the ``canonical''
 matrix from \eqref{DN}.

 Similar to \eqref{repl1}--\eqref{replN2},
 replacement rules can be applied to polynomials~$C_N$.

In the above notation $ b(\mathfrak{A})$ from the solution of
system  \eqref{cPis0} can be expressed as
\begin{equation}\label{bR1}
  b(\mathfrak{A})=Q(\A)
\end{equation}
where
\begin{equation}\label{R1}
  Q(\{a_0,\dots,a_{N-1}\})=\frac{D_N\big(x_{j,k}\rightarrow
  \eta^{\langle k \rangle }(a_j)\big)}
  {D_N\big(x_{j,0}\rightarrow 1,x_{j,k}\rightarrow
  \eta^{\langle k \rangle }(a_j)\big)}.
\end{equation}

For the sets from
Fig.~\ref{gc} and Tables \ref{manygrids}--\ref{manycircles}
we have
\begin{equation}\label{R1Ais1}
  Q(\mathfrak{A})\approx  1.
\end{equation}
Of course, this cannot be true for
an arbitrary set $\mathfrak{A}$.
To see this, let, for example, $\mathfrak{A}_N(a)=\{a,2a,\dots,(N-1)a\}$.
Then $Q(\mathfrak{A}_N(a))$ is a meromorphic
function,  different from constant (otherwise we would have a differential equation for $\eta(s)$).
Hence, by Piccard's theorems, all
(except, possibly, one)
complex numbers are  values of~$Q(\mathfrak{A}_N(a))$.

It is not clear in what terms could we describe the
class of sets $\mathfrak{A}$ for which
\begin{equation}
  \vert Q(\mathfrak{A})-  1\vert < \varepsilon
\end{equation}
for a given positive~$\varepsilon$.
A more
achievable goal could be to find bounds for cases of sets of particular structure, for example,
to establish inequalities of the form
\begin{equation}
  \vert Q(\mathfrak{A}_\mathrm{G}(a,\delta_1,\delta_2,N_1,N_2))
  -1 \vert < B_\mathrm{G}(a,\delta_1,\delta_2,N_1,N_2)
\end{equation}
and
\begin{equation}
   \vert Q(\mathfrak{A}_\mathrm{C}(c,r,N)
   -1\vert<B_\mathrm{C}(c,r,N)
\end{equation}
with certain functions  $B_\mathrm{G}(a,\delta_1,\delta_2,N_1,N_2)$ and
$B_\mathrm{C}(c,r,N)$.

\section{Two generalization}\label{sect1G}

Approximate inequality \eqref{R1Ais1}
can be generalized in several directions.

\subsection{Characteristic
polynomials}\label{charpol}
The
determinant of a matrix is just the constant term
of its characteristic
polynomial;
in our notation
\begin{equation}\label{DisC0}
D_N=C_N(\lambda\rightarrow 0).
\end{equation}
In particular, the numerator and the denominator in
\eqref{R1} are the constant terms of the following
 polynomials:
\begin{equation}\label{V1}
  {C_N\big( x_{j,k}\rightarrow
  \eta^{\langle k \rangle }(a_j)\big)}
  =\\V_{\Nc N}(\A)\lambda^{N}+\dots+V_{\Nc n}(\A)\lambda^{n}
  +\dots+V_{\Nc 0}(\A),
 \end{equation}
 \begin{multline}\label{W1}
   {C_N\big(  x_{j,0}\rightarrow 1,\ x_{j,k}\rightarrow
  \eta^{\langle k \rangle }(a_j)\big)}=\\
  W_{\Nc N}(\A)\lambda^{N}+\dots+
  W_{\Nc n}(\A)\lambda^{n}+\dots+W_{\Nc 0}(\A).
\end{multline}

\begin{table}[t]
\centering
\input gridcharpolnew.tex
\caption{Generalization of  \eqref{R1Ais1}
for grid  $\A_1$ from Fig.~\ref{gc}.}
\label{gridcharpol}
\end{table}

Let us consider ratios
\begin{equation}\label{Qnis1}
  Q_{n}(\mathfrak{A})=
  \frac{V_{\Nc n}(\A)}{W_{\Nc n}(\A)}, \quad n=0,\dots N.
\end{equation}
Clearly,
$ Q_{0}(\mathfrak{A})= Q(\mathfrak{A})$,
and we have seen that this ratio is close to~$1$
for  the sets from Tables~\ref{manygrids}--\ref{manycircles} and for sets $\mathfrak{A}_1 $ and $\mathfrak{A}_2 $ from Fig.\ \ref{gc}.
Table~\ref{gridcharpol} demonstrates that
 other ratios    $Q_{n}(\mathfrak{A}_1) $
are also not far from~$1$ when $n$ is small.

\subsection{Missing derivatives}\label{misder}

In our search for (an impossible) relation
\eqref{Pis0} we could in advance exclude
certain derivatives. More formally,
let  $\D=d_0,d_1,\dots$
be a strictly increasing sequences
of non-negative integers and
\begin{equation}\label{d00}
  d_0=0.
\end{equation}
We consider now the following generalization
of \eqref{cPis0}:
\begin{equation}\label{cPis0gen}
  \eta(a_j)+c_{d_1}(\A,\D)\eta^{\langle d_1 \rangle}(a_j)+
  \dots+c_{d_{N-1}}(\A,\D)\eta^{\langle d_{N-1} \rangle}(a_j)=\\b(\A,\D).
\end{equation}
In the solution of this system (with $j=0,\dots,N-1$)
\begin{equation}\label{bR1gen}
  b(\mathfrak{A},\D)=
  Q(\mathfrak{A},\D)=\frac{D_N\big(x_{j,k}\rightarrow
  \eta^{\langle d_k \rangle }(a_j)\big)}
  {D_N\big(x_{j,0}\rightarrow 1,x_{j,k}\rightarrow
  \eta^{\langle d_k \rangle }(a_j)\big)}.
\end{equation}

A sequence $\D$ can be defined by
the  set $\M_\D$ of non-negative integers missing in~$\D$.
Table \ref{circle_drop} shows how dropping more
and more derivatives increases the distance
between $Q(\mathfrak{A},\D)$ and $1$ for
a particular set $\A$ and a selection of
sets $\M_\D$ of missing derivatives.

\begin{table}[t]
\centering
\input circle_drop_new.tex
\caption{The distance between
$Q(\A,\D)$ and $1$
for $\A=\A_{\mathrm{C}}(0.6+14\mathrm{i},10,20)$.
}
\label{circle_drop}
\end{table}

\section{Dropping $\eta(a)$}\label{dropeta}

Searching for
a polynomial
\begin{equation}\label{polP0}
  P(y_0,y_1,\dots,y_{N-1})=
  c_0y_0+c_1y_1+\dots+c_{N-1}y_{N-1}-b
\end{equation}
such that
\begin{equation}\label{P0is0}
  P\big(\eta(a),\eta'(a),\dots,
  \eta^{\langle N-1 \rangle}(a)\big)=
  0
\end{equation}
we should exclude the degenerate case
of polynomial identically equal to $0$.
In \eqref{polP} this was achieved
by fixing~$c_0=1$ but now we remove this constraint.
Let $l$ be the least index such that $c_l\ne0$;
without loss of generality we put $c_l=1$.

More formally, for a given  set $\A=\{a_0,\dots,a_{N-1}\}$
and non-negative integer $l$
we generalize  \eqref{cPis0} to
linear system
 \begin{equation}\label{2cPis0}
  \eta^{\langle l \rangle}(a_j)+c_{l+1}(\A)\eta^{\langle l+1 \rangle}(a_j)+
  \dots\\+c_{l+N-1}(\A)\eta^{\langle l+N-1 \rangle}(a_j)=b(\A)
\end{equation}
with $j=0,\dots,N-1$.
We find that the case $l\ge 1$
radically differs from the earlier
considered case $l=0$. Namely,
$b(\A)$  has the tendency to be
close to~$0$.
In our notation this can be written as
 \begin{eqnarray}\label{appr2}
  {D_{N}\big(x_{j,k}\rightarrow
  \eta^{\langle l+k \rangle }(a_j)\big)}
  \approx 0.
 \end{eqnarray}
However, this approximate equality
  is
 uninformative.
Namely,
polynomial
$ D_N$
 has no constant term,
so the value of  the left-hand side of
\eqref{appr2} can be small
just due to the small values of the eta
derivatives.

Nevertheless, approximate equality
\eqref{appr2} can provide non-trivial
information about the relationship
between the derivatives of the eta function.
  To demonstrate this,
  we  select two integers $m$ and $n$ such that
 $0\le m<N$ and $0\le n<N$  and define $y$
 as  the solution of linear equation
 \begin{eqnarray}\label{appr2eq}
  {D_{N}\big(x_{m,n}\rightarrow y,\ x_{j,k}\rightarrow
  \eta^{\langle l+k \rangle }(a_j)\big)}
  = 0.
 \end{eqnarray}
   In our notation,
 \begin{multline}\label{R2}
   y=
   R({l,m,n},\A)=\\
   \frac{D_{N}\big(x_{m,n}\rightarrow 0,\ x_{j,k}\rightarrow
  \eta^{\langle l+k \rangle }(a_j)\big)}
  {D_{N}\big(x_{m,n}\rightarrow -1,\ x_{m,k}\rightarrow
  0,\ x_{j,n}\rightarrow
  0,\ x_{j,k}\rightarrow
  \eta^{\langle l+k \rangle }(a_j)\big)}.
 \end{multline}
According to \eqref{appr2}--\eqref{appr2eq},
 we can expect that
 $y\approx  \eta^{\langle l+n \rangle }(a_m)$.

\begin{example}\normalfont
   Let $\A=\{a_0,\dots,a_{11}\}$ be either of the two sets from Fig.~\ref{gc}. We have:
  for $l=1,\dots,11$, $m=0,\dots,11$,  $n=0,\dots,11$
\begin{equation}\label{rel2example}
  \left\vert\frac{ R({l,m,n},\A)}{\eta^{\langle l+n \rangle }(a_m)}- 1\right\vert<10^{-7}.
\end{equation}
\end{example}

\begin{remark}\normalfont
 As in Subsection \ref{charpol}
  we can consider
      the characteristic polynomials
 \begin{equation}\label{eig2nbis}
   {C_{N}\big(x_{m,n}\rightarrow 0,\ x_{j,k}\rightarrow
  \eta^{\langle l+k \rangle }(a_j)\big)}
 \end{equation}
and
 \begin{equation}\label{eig2d2bis}
  {C_{N}\big(x_{m,n}\rightarrow -1,\ x_{m,k}\rightarrow
  0,\ x_{j,n}\rightarrow
  0,\ x_{j,k}\rightarrow
  \eta^{\langle l+k \rangle }(a_j)\big)}
 \end{equation}
 and approximate
 $\eta^{\langle l+n \rangle}(a)$ by ratios of
   lower coefficients of
   these polynomials.
\end{remark}

\begin{remark}\normalfont
 As in Subsection \ref{misder}
  we could drop not only several initial
 but also some higher derivatives.
\end{remark}

\begin{remark}\normalfont
Similar to \eqref{eig2nbis}--\eqref{eig2d2bis}
we can transform \eqref{R1}--\eqref{R1Ais1}
into linear equation
\begin{equation}\label{R1y}
{D_N\big(x_{m,n}\rightarrow y, x_{j,k}\rightarrow
  \eta^{\langle k \rangle }(a_j)\big)}=
  {D_N\big(x_{m,n}\rightarrow y,
  x_{j,0}\rightarrow 1,x_{j,k}\rightarrow
  \eta^{\langle k \rangle }(a_j)\big)}
\end{equation}
and  expect that its solution
is close to  $\eta^{\langle n \rangle }(a_m)$.

\end{remark}

\part{Conjectures}\label{part2}
\section{First limiting case}\label{1lim}

\begin{table}[t]
\centering
\scriptsize
\input R3_new.tex
\normalsize
\caption{Numerical data supporting
Conjecture  \ref{conj1}.}
\label{R3table}
\end{table}

Approximate equality \eqref{R1Ais1} can be  fulfilled with high
accuracy when points $a_0,\dots,a_{N-1}$ are spaced apart
(see the cases of  large
$\delta_1$ or $\delta_2$ in
Table~\ref{manygrids} and \cite[Table 36]{preprint},
and also
the cases of  large
$r$ in
Table \ref{manycircles}  and  \cite[Table 37]{preprint}). But  the accuracy of
\eqref{R1Ais1} can also  be high
when the points $a_0,\dots,a_{N-1}$
are close to each other (see the cases when $r$ or both
$\delta_1$ and $\delta_2$
 are small
in the same tables).
Thus we can consider the limit of
$Q(\A)$ when some of the points from set $\A$ approach
the same limiting values.

In this paper we restrict ourselves to the extreme case
when all the points tend to the same number.
Let $a_j=a+j\varepsilon$ where $\varepsilon$ tends to~$0$.
We have:
 \begin{equation}\label{lim1num}
  D_N\big(x_{j,k}\rightarrow
  \eta^{\langle k \rangle }(a_j)\big)=
  D_N(x_{j,k}\rightarrow
  \eta^{\langle j+k \rangle }(a)\big)\varepsilon^{\frac{N(N-1)}{2}}+
  O(\varepsilon^{\frac{N(N-1)}{2}+1}),
\end{equation}
 \begin{multline}\label{lim1den}
  D_N\big(x_{j,0}\rightarrow 1,x_{j,k}\rightarrow
  \eta^{\langle k \rangle }(a_j)\big)=\\
  {D_{N-1}\big(x_{j,k}\rightarrow
  \eta^{\langle j+k+2 \rangle }(a)\big)\varepsilon^{\frac{N(N-1)}{2}}}+
  O(\varepsilon^{\frac{N(N-1)}{2}+1}).
\end{multline}

Let
\begin{eqnarray}\label{R3}
  S_N(a)&=&\frac
  {D_N\big(x_{j,k}\rightarrow
  \eta^{\langle j+k \rangle }(a)\big)}
  {D_{N-1}\big(x_{j,k}\rightarrow
  \eta^{\langle j+k+2 \rangle }(a)\big)}
\end{eqnarray}
From \eqref{lim1num} and \eqref{lim1den}
\begin{eqnarray}\label{R3aeps}
  S_N(a)
 &=& Q\big(\{a_0,\dots,
  a_{N-1}\}\big)+O(\varepsilon ),
\end{eqnarray}
and according to \eqref{R1Ais1}  we can expect that $S_N(a)$ is
   close to~$1$.

\begin{conjecture}\label{conj1}
Except for countably many values of $a$,
\begin{equation}\label{R3is1}
   S_N(a)
   \stackrel[N\rightarrow \infty]{\longrightarrow}{} 1
\end{equation}
where  $S$
is defined by
\eqref{R3}.
\end{conjecture}

Numerical data in Table \ref{R3table}  support this conjecture
(for other relevant data see \cite[Table 1]{preprint}).

\begin{remark}\normalfont
  Similar to \eqref{appr2eq}, we can consider linear equation
  in unknown~$y$
  \begin{eqnarray}\label{R3eq}
    {D_N\big(x_{0,0}\rightarrow y,\ x_{j,k}\rightarrow
  \eta^{\langle j+k \rangle }(a)\big)}
  ={D_{N-1}\big(x_{j,k}\rightarrow
  \eta^{\langle j+k+2 \rangle }(a)\big)}.
\end{eqnarray}
Conjecture \ref{conj1} suggests that
the solution of this equation should be close to~$\eta(a)$.
Other ways to calculate $\eta(a)$ via derivatives of
  $\eta(s)$ at $s=a$ were proposed
  in \cite[Section 6]{NO0preprint} (see also \cite[Section 8)]{Sirius}).

\end{remark}

\section{Second  limiting case}\label{2lim}

\begin{table}[t]
\centering
\scriptsize
\input conj3new.tex
\normalsize
\caption{Numerical data supporting
Conjecture \ref{conj3}; $N=50$, $a=0.4+17\myi$.}
\label{conj3data}
\end{table}

Now we consider limiting behavior of
\eqref{appr2}.
Similar to \eqref{lim1num} we have:
 \begin{multline}\label{lim2num}
  D_N\big(x_{j,k}\rightarrow
  \eta^{\langle l+k \rangle }(a+j\varepsilon)\big)=\\
  D_N\big(x_{j,k}\rightarrow
  \eta^{\langle l+j+k \rangle }(a)\big)\varepsilon^{\frac{N(N-1)}{2}}+
  O(\varepsilon^{\frac{N(N-1)}{2}+1}).
\end{multline}
We can make a guess that
 \begin{equation}\label{lim2lguess}
  D_N\big(x_{j,k}\rightarrow
  \eta^{\langle l+j+k \rangle }(a)\big)\approx 0.
\end{equation}

Like \eqref{appr2}, by itself this is an uninformative approximate relation
among the derivatives of the eta function.
To show its non-triviality we
can consider counterparts of
\eqref{appr2eq},
namely, equations
 \begin{equation}\label{lim2leq}
  D_N\big(j+k=m\Rightarrow x_{j,k}\rightarrow y,
  x_{j,k}\rightarrow
  \eta^{\langle l+j+k \rangle }(a)\big)= 0
\end{equation}
for $m=0,\dots,2N-2$.
However, there is an essential distinction:
\eqref{appr2} is a
linear equation in $y$
but \eqref{lim2leq} is
 linear
 in  two extreme cases only,
 namely, when $m=0$ and when $m=2N-2$.
For these values of $m$ the solution
 is given by
\begin{equation}\label{R4}
   y=
   \frac{D_{N}\big(x_{m',m'}\rightarrow 0,\ x_{j,k}\rightarrow
  \eta^{\langle l+j+k \rangle }(a)\big)}
  {D_{N}\big(x_{m',m'}\rightarrow -1,\ x_{m',k}\rightarrow
  0,\ x_{j,m'}\rightarrow
  0,\ x_{j,k}\rightarrow
  \eta^{\langle l+j+k \rangle }(a)\big)}
 \end{equation}
where $m'=m/2$.

 In general case, the degree of
 equation \eqref{lim2leq} is equal to
 \begin{equation}
  d=\min(m+1,2N-m-1).
 \end{equation}
  Let
  \begin{align}
\hspace{-3mm}E_{l,m,N}(a,y)&={D_N\big(j+k=m\Rightarrow x_{j,k
}\rightarrow y,\
  x_{j,k}\rightarrow
  \eta^{\langle l+j+k \rangle }(a)\big)}\nonumber\\
  \label{EE}
  &=
   E_{l,m,N,d}(a)y^d+\dots+E_{l,m,N,n}(a)y^n
   +\dots+E_{l,m,N,0}(a).
\end{align}

Assuming  \ref{lim2lguess}, we can expect that
\emph{at least one} of the roots of equation~\eqref{lim2leq}
is close to $\eta^{\langle l+m \rangle}(a)$.
Calculations show that when $d$ is relatively
small to $N$, \emph{all}  roots of
equation \eqref{lim2leq}
are close to~$\eta^{\langle l+m \rangle}(a)$.

Here we consider only the case when $d$ small
due to the small value of~$m$ (for the case when
$m$ is close to $2N-2$ see
 \cite[Subsection 3.4]{preprint}).

\begin{conjecture} \label{conj2}
Except for countably many values of $a$, for every
positive integer $l$ and
non-negative integer $m$   the maximal distance between
   the roots of polynomial $  E_{l,m,N}(a,y)$
(defined by \eqref{EE}) and $\eta^{\langle l+m \rangle}(a)$
tends to zero as~\mbox{$N\rightarrow\infty$.}
\end{conjecture}

This conjecture does not have
the ``canonical'' form~\eqref{VWlim1},
and we do not present here numerical examples
which would
support it (for such data  see
\cite[Tables 11-16]{preprint}. Instead of this, we deduce
several
corollaries of the form~\eqref{VWlim1}.
Taken together, they
are essentially equivalent to
Conjecture~\ref{conj2}.

Namely, Conjecture~\ref{conj2}
implies that if $m$ is small (relative to $N$) then
 polynomial $E_{l,m,N}(a,y)$
should have a special structure, namely,
\begin{equation}
   E_{l,m,N}(a,y)\approx A(y-\eta^{\langle l+m \rangle}(a))^d
\end{equation}
for some non-zero constant~$A$. This opens
several ways to calculate
(approximation to)~$\eta^{\langle m \rangle}(a)$
via the values of other derivatives.

\begin{conjecture} \label{conj3}
Except for countably many values of $a$, for every
positive integer $l$, non-negative integer $m$ and
an  $n$ such that $0<n\le m+1$
\begin{equation}\label{etaklim}
  \frac{n-m-2}{n}\cdot\frac{E_{l,m,N,n-1}(a)}{E_{l,m,N,n}(a)}
  \stackrel[N\rightarrow \infty]{\longrightarrow}{}\eta^{\langle l+m \rangle}(a)
\end{equation}
where $E_{l,m,N,n}$
is defined by
\eqref{EE}.
\end{conjecture}

For numerical data supporting Conjecture \ref{conj3}
see Table \ref{conj3data} and \cite[Tables 17--19]{preprint}.

\begin{remark}\normalfont
Assuming Conjecture  \ref{conj1},
we can consider equations
\begin{multline}
{D_N\big(j+k=m\Rightarrow x_{j,k}\rightarrow y,\
  x_{j,k}\rightarrow
  \eta^{\langle j+k \rangle }(a)\big)}=\\
  {D_{N-1}\big(j+k+2=m\Rightarrow x_{j,k}\rightarrow y,\ x_{j,k}\rightarrow
  \eta^{\langle j+k+2 \rangle }(a)\big)}\label{E1}
  \end{multline}
and state for them analogs of Conjectures
\ref{conj2} and \ref{conj3}  (for numerical data
see \cite[Tables 7-10]{preprint}.

\end{remark}

\section{Modifications}\label{vary}

Besides \eqref{R3is1} and
\eqref{etaklim} (and corresponding
 generalizations indicated in the remarks),
one can state similar conjectures
of the form \eqref{VWlim1} for
a great number of  other
ratios of polynomial.
Here we indicate a few directions of
possible modifications.

\subsection{Taylor coefficients}

\begin{table}[t]
\centering
\scriptsize
\input R3fact_new.tex
\normalsize
\caption{Numerical data supporting
Conjecture  \ref{conj1fact}.}
\label{R3facttable}
\end{table}

\begin{table}[t]
\centering
\scriptsize
\input Rfact3.tex
\normalsize
\caption{Numerical data supporting
Conjecture  \ref{conjfact1}.}
\label{factR3table}
\end{table}

Let
\begin{equation}\label{DNfact}
  D_{!\,N}=D_N(x_{j,k}\rightarrow x_{j,k}/k!).
\end{equation}
By analogy with
\eqref{R3} we define
\begin{eqnarray}\label{R3fact}
  S_{!\,N}(a)&=&\frac
  {D_{!\,N}\big(x_{j,k}\rightarrow
  \eta^{\langle j+k \rangle }(a)\big)}
  {D_{!\,N-1}\big(x_{j,k}\rightarrow
  \eta^{\langle j+k+2 \rangle }(a)\big)}
\end{eqnarray}
and state a counterpart of Conjecture \ref{conj1}.

\begin{conjecture}\label{conj1fact}
Except for countably many values of $a$,
\begin{equation}\label{conj4}
   S_{!\,N}(a)
   \stackrel[N\rightarrow \infty]{\longrightarrow}{} 1.
\end{equation}
\end{conjecture}

Comparing values of $S_{!\,N}(a)$ in Table~\ref{R3facttable}
(see also \cite[Tables 27]{preprint})
with values of  $S_{N}(a)$ in Table~\ref{R3table}
for the same  $N$ and $a$, we
should observe that convergence
in \eqref{conj4} is much slower than in~\eqref{R3is1}.

\subsection{``Anty-Taylor'' coefficients}

Similar to \eqref{DNfact} and \eqref{R3fact}, let
\begin{equation}\label{factDN}
  D^{!}_{N}=D_N(x_{j,k}\rightarrow k! x_{j,k})
\end{equation}
and
\begin{eqnarray}\label{Rfact3}
  S^{!}_{N}(a)&=&\frac
  {D^{!}_{N}\big(x_{j,k}\rightarrow
  \eta^{\langle j+k \rangle }(a)\big)}
  {D^{!}_{N-1}\big(x_{j,k}\rightarrow
  \eta^{\langle j+k+2 \rangle }(a)\big)}.
\end{eqnarray}

\begin{conjecture}\label{conjfact1}
Except for countably many values of $a$,
\begin{equation}\label{Rfact3is1}
   S^{!}_{N}(a)
   \stackrel[N\rightarrow \infty]{\longrightarrow}{} 1.
\end{equation}
\end{conjecture}

Comparing values of $S^{!}_{N}(a)$ in Table~\ref{factR3table}
(see also \cite[Tables 29]{preprint})
with values of $S_{!\, N}(a)$ in Table~\ref{R3facttable}
for the same  $N$ and $a$, we
should apprehend that convergence
in \eqref{Rfact3is1} is even slower than in~\eqref{conj4}.

\subsection{Minors}

\begin{table}
\centering
\scriptsize
\input 1stmissing.tex
\normalsize
\caption{Numerical data supporting
Conjecture  \ref{conj6}; $a=0.4+16\myi$.}
\label{dropdrop}
\end{table}

In this subsection we extend Conjecture \ref{conj1}.
Let
\begin{equation}\label{dseq}
 \D=d_0,d_1,\dots
\end{equation}
 and
\begin{equation}\label{fseq}
  \F=f_0,f_1\dots
\end{equation}
be two strictly increasing sequences
of non-negative integers.
We generalize~\eqref{DN} to
\begin{equation}\label{DNDF}
  D_{N,\D,\F}=\begin{vmatrix}
    x_{f_0,d_0}&\dots&x_{f_0,d_{N-1}}\\
    \vdots&\ddots&\vdots\\
     x_{f_{N-1},d_0}&\dots&x_{f_{N-1},d_{N-1}}
  \end{vmatrix}.
\end{equation}

As in
Subsection \ref{misder},
$\M_\D$ and $\M_\F$ denote the sets of
non-negative integers missing in
$\D$ and $\F$ respectively.

\begin{conjecture}\label{conj6}
Let $\D_1=d_1,d_2,\dots$ and  $\F_1=f_1,f_2\dots$
be the tailors of sequences \eqref{dseq} and~\eqref{fseq}.
If $d_0=f_0=0$ and sets $\M_\D$ and $\M_\F$ are finite,
then except for countably many values of $a$,
\begin{eqnarray}\label{R3DF}
  S_{N,\D,\F}(a)&=&\frac
  {D_{N,\D,\F}\big(x_{j,k}\rightarrow
  \eta^{\langle j+k \rangle }(a)\big)}
  {D_{N-1,\D_1,\F_1}\big(x_{j,k}\rightarrow
  \eta^{\langle j+k+2 \rangle }(a)\big)}
  \stackrel[N\rightarrow \infty]{\longrightarrow}{} 1.
\end{eqnarray}
\end{conjecture}

Table \ref{dropdrop} presents data supporting Conjecture \ref{conj6}
(see also \cite[Table 32]{preprint}).

\begin{remark}\normalfont
Very probably, the requirement of the finiteness of
sets $\M_\D$ and $\M_\F$  is excessive.
Most likely, it
be replaced  by certain restrictions on the growth
of sequences \eqref{dseq} and \eqref{fseq}
(cf. \cite[Tables 33--34]{preprint}).

\end{remark}

\subsection{Generalized characteristic polynomials  }

\begin{table}[t]
\centering
\scriptsize
\input gen_char_poly.tex
\normalsize
\caption{Approximation of $\eta^{\langle l \rangle}(a)$
by ratios of the    initial
coefficients of polynomials
\eqref{V1bis} and \eqref{W1bis}; $l=2$, $a=0.7+19\myi$.}
\label{genchar}
\end{table}

Let
\begin{equation}\label{alpha}
  \G=\langle  g_0,\dots,g_{N-1}\rangle
\end{equation}
 be
a permutation of numbers
$0,\dots,N-1$. We generalize  \eqref{C}
to
\begin{equation}\label{CG}
  C_{\Nc\G}=D_N(j=g_k\Rightarrow x_{j,k}\rightarrow x_{j,k}-\lambda).
\end{equation}
Similar to \eqref{DisC0},
$D_N$ is the constant term of  $C_{\Nc\G}$
up to the sign (depending on the parity of $\G$).
Respectively, for $m=0$  ratio \eqref{R4} is equal to the ratio
of the constant terms of polynomials
\begin{multline}\label{V1bis}
  C_{\Nc\G}\big(x_{0,0}\rightarrow 0,\ x_{j,k}\rightarrow
  \eta^{\langle l+j+k \rangle }(a)\big)
  =\\V_{l,N,\G}(a)\lambda^{N}+\dots+V_{l,n,\G}(a)\lambda^{n}
  +\dots+V_{l,0,\G}(a),
 \end{multline}
 \begin{multline}\label{W1bis}
   C_{\Nc\G}\big(x_{0,0}\rightarrow -1,\ x_{0,k}\rightarrow
  0,\ x_{j,0}\rightarrow
  0,\ x_{j,k}\rightarrow
  \eta^{\langle l+j+k \rangle }(a)\big)
  =\\W_{l,N,\G}(a)\lambda^{N}+\dots+W_{l,n,\G}(a)\lambda^{n}
  +\dots+W_{l,0,\G}(a).
\end{multline}

For relatively prime positive integers $q$ and $r$ let
$\alpha_{q,r}$ denote such permutation \eqref{alpha} that
\begin{equation}
  g_k\equiv kr \pmod{q}.
\end{equation}
Similar to Table \ref{gridcharpol},
we observe in Table \ref{genchar}
that for small $n$  ratios
$V_{l,n,\G}(a)/W_{l,n,\G}(a)$ are close to
$\eta^{\langle l \rangle}(a)$.
For other relevant data see
\cite[Subsection 4.6]{preprint}.

\section*{Concluding remark}\label{zetasection}

All our numerical data were for the eta function.
Very likely, other functions defined
by Dirichlet series
have similar properties.

Of course, the most interesting case is the classic
Riemann zeta function. Calculations show that
for this function we can observe phenomena
similar to those discussed in Part \ref{part1}.
However, their accuracy can be rather low when
some elements of $\A$ are close to the
pole of the zeta function.
Calculations performed so far do not allow us to
put forward  counterparts of
Conjectures~\ref{conj1}--\ref{conj6} (for the relevant numerical data see
\cite[Tables 39--40]{preprint}).

\section*{Acknowledgments}

Numerical data for Tables 1--11 were computed
by programs written in {\sc \mbox{Julia}}
programming language~\cite{JULIA}
with package {\sc Nemo} \cite{Nemo}.

\printbibliography

Author

\

Yuri Matiyasevich

Steklov Institute of Mathematics at St.Petersburg

27 Fontanka, St.Petersburg, 191023 Russia

yumat [at] pdmi [dot] ras [dot] ru

\url{http://logic.pdmi.ras.ru/~yumat}

\end{document}

%% file: gccoeffs_2.tex
\begin{tabular}{||c||r|r||}
  \hhline{|t:=:t:==:t|}
    \multicolumn{1}{||c||}{}&\multicolumn{1}{c|}{$\mathfrak{A}=\mathfrak{A}_1$}&\multicolumn{1}{c||}{$\mathfrak{A}=\mathfrak{A}_2$}
  \\ 
  \hhline{|:=::==:|}
    $c_{1}(\mathfrak{A})$
  &
    $\ \ 7.30910\ldots\ \ -0.26348\ldots\mathrm{i}$
  &
    $\ \ 7.35147\ldots\ \ -0.23197\ldots\mathrm{i}$
  \\ 
    $c_{2}(\frak{A})$
  &
    $\ 23.85173\ldots\ \ -1.81924\ldots\mathrm{i}$
  &
    $\ 24.15410\ldots\ \ -1.61238\ldots\mathrm{i}$
  \\ 
    $c_{3}(\frak{A})$
  &
    $\ 45.96860\ldots\ \ -5.55318\ldots\mathrm{i}$
  &
    $\ 46.92441\ldots\ \ -4.95964\ldots\mathrm{i}$
  \\ 
    $c_{4}(\frak{A})$
  &
    $\ 58.22426\ldots\ \ -9.89546\ldots\mathrm{i}$
  &
    $\ 59.99059\ldots\ \ -8.91557\ldots\mathrm{i}$
  \\ 
    $c_{5}(\frak{A})$
  &
    $\ 50.94334\ldots\ -11.42336\ldots\mathrm{i}$
  &
    $\ 53.05986\ldots\ -10.39502\ldots\mathrm{i}$
  \\ 
    $c_{6}(\frak{A})$
  &
    $\ 31.43810\ldots\ \ -8.94113\ldots\mathrm{i}$
  &
    $\ 33.15830\ldots\ \ -8.22831\ldots\mathrm{i}$
  \\ 
    $c_{7}(\frak{A})$
  &
    $\ 13.68711\ldots\ \ -4.81140\ldots\mathrm{i}$
  &
    $\ 14.64820\ldots\ \ -4.48445\ldots\mathrm{i}$
  \\ 
    $c_{8}(\frak{A})$
  &
    $\ \ 4.11913\ldots\ \ -1.75935\ldots\mathrm{i}$
  &
    $\ \ 4.48375\ldots\ \ -1.66349\ldots\mathrm{i}$
  \\ 
    $c_{9}(\frak{A})$
  &
    $\ \ 0.81560\ldots\ \ -0.41866\ldots\mathrm{i}$
  &
    $\ \ 0.90549\ldots\ \ -0.40230\ldots\mathrm{i}$
  \\ 
    $c_{10}(\frak{A})$
  &
    $\ \ 0.09551\ldots\ \ -0.05857\ldots\mathrm{i}$
  &
    $\ \ 0.10851\ldots\ \ -0.05731\ldots\mathrm{i}$
  \\ 
    $c_{11}(\frak{A})$
  &
    $\ \ 0.00500\ldots\ \ -0.00365\ldots\mathrm{i}$
  &
    $\ \ 0.00583\ldots\ \ -0.00365\ldots\mathrm{i}$
  \\ 
  \hhline{|b:=:b:==:b|}
\end{tabular}

%% file: multigrids.tex
\begin{tabular}{||r|r|r|r|r||r|l||}
  \hhline{|t:=====:t:==:t|}
    \multicolumn{1}{||c|}{$a$}&\multicolumn{1}{c|}{$\delta_1$}&\multicolumn{1}{c|}{$\delta_2$}&\multicolumn{1}{c|}{$N_1$}&\multicolumn{1}{c||}{$N_2$}&\multicolumn{1}{c|}{$N$}&\multicolumn{1}{c||}{\vpad{$\vert b(\frak{A})-1\vert$}}
  \\ 
  \hhline{|:=====::==:|}
    $-200+15\mathrm{i}$
  &
    $50$
  &
    $50\mathrm{i}$
  &
    $5$
  &
    $5$
  &
    $36$
  &
    $\ 1.0523\dotsc\cdot 10^{-57}$
  \\ 
    $-5+15\mathrm{i}$
  &
    $5$
  &
    $5\mathrm{i}$
  &
    $5$
  &
    $5$
  &
    $36$
  &
    $\ 5.7251\dotsc\cdot 10^{-64}$
  \\ 
    $30\mathrm{i}$
  &
    $0.0001$
  &
    $0.0001\mathrm{i}$
  &
    $5$
  &
    $5$
  &
    $36$
  &
    $\ 1.1424\dotsc\cdot 10^{-53}$
  \\ 
    $0.25+30\mathrm{i}$
  &
    $0.0001$
  &
    $0$
  &
    $59$
  &
    $0$
  &
    $60$
  &
    $\ 5.3809\dotsc\cdot 10^{-98}$
  \\ 
    $0.5+30\mathrm{i}$
  &
    $0$
  &
    $3\mathrm{i}$
  &
    $0$
  &
    $59$
  &
    $60$
  &
    $\ 2.0739\dotsc\cdot 10^{-56}$
  \\ 
    $-200+10000\mathrm{i}$
  &
    $50$
  &
    $50\mathrm{i}$
  &
    $5$
  &
    $5$
  &
    $36$
  &
    $\ 3.8883\dotsc\cdot 10^{-51}$
  \\ 
  \hhline{|b:=====:b:==:b|}
\end{tabular}

%% file: multicircles.tex
\begin{tabular}{||r|r|r||l||}
  \hhline{|t:===:t:=:t|}
    \multicolumn{1}{||c|}{$c$}&\multicolumn{1}{c|}{$r$}&\multicolumn{1}{c||}{$N$}&\multicolumn{1}{c||}{\vpad{$\vert b(\mathfrak{A}_\mathrm{C}(a,
r,N))-1\vert$}}
  \\ 
  \hhline{|:===::=:|}
    $10\mathrm{i}$
  &
    $10^{-5}$
  &
    $24$
  &
    $\ 1.2861\dotsc\cdot 10^{-32}$
  \\ 
    $-1+20\mathrm{i}$
  &
    $10^{-3}$
  &
    $36$
  &
    $\ 6.2669\dotsc\cdot 10^{-52}$
  \\ 
    $0.25+30\mathrm{i}$
  &
    $10^{-1}$
  &
    $60$
  &
    $\ 5.4442\dotsc\cdot 10^{-98}$
  \\ 
    $0.6+40\mathrm{i}$
  &
    $1$
  &
    $72$
  &
    $\ 1.1660\dotsc\cdot 10^{-121}$
  \\ 
    $1+50\mathrm{i}$
  &
    $10$
  &
    $80$
  &
    $\ 4.6580\dotsc\cdot 10^{-138}$
  \\ 
    $2+60\mathrm{i}$
  &
    $10^{2}$
  &
    $96$
  &
    $\ 4.5686\dotsc\cdot 10^{-225}$
  \\ 
    $70\mathrm{i}$
  &
    $10^{3}$
  &
    $120$
  &
    $\ 3.9966\dotsc\cdot 10^{-1623}$
  \\ 
  \hhline{|b:===:b:=:b|}
\end{tabular}

%% file: gridcharpolnew.tex
\begin{tabular}{||r||rlrl|rl||}
  \hhline{|t:=:t:======:t|}
    $n$&\multicolumn{4}{c|}{$V_{n}(\A_1)$}&\multicolumn{2}{c||}{\vpad{$\vert Q_{n}(\A_1)-1\vert$}}
  \\ 
  \hhline{|:=::======:|}
    $0$
  &
    $-2.8932\dotsc\cdot$
  &
    $\hspace{-4mm} 10^{-48}$
  &
    $\hspace{-4mm}-5.1873\dotsc\cdot$
  &
    $\hspace{-4mm} 10^{-49}\myi$
  &
    $\hspace{9mm} 1.1085\dotsc\cdot$
  &
    $\hspace{-4mm} 10^{-14}$
  \\ 
    $1$
  &
    $1.1853\dotsc\cdot$
  &
    $\hspace{-4mm} 10^{-33}$
  &
    $\hspace{-4mm}+2.1097\dotsc\cdot$
  &
    $\hspace{-4mm} 10^{-34}\myi$
  &
    $\hspace{9mm} 8.6263\dotsc\cdot$
  &
    $\hspace{-4mm} 10^{-13}$
  \\ 
    $2$
  &
    $-2.2690\dotsc\cdot$
  &
    $\hspace{-4mm} 10^{-22}$
  &
    $\hspace{-4mm}-1.6178\dotsc\cdot$
  &
    $\hspace{-4mm} 10^{-22}\myi$
  &
    $\hspace{9mm} 5.1025\dotsc\cdot$
  &
    $\hspace{-4mm} 10^{-11}$
  \\ 
    $3$
  &
    $4.0262\dotsc\cdot$
  &
    $\hspace{-4mm} 10^{-14}$
  &
    $\hspace{-4mm}+8.4755\dotsc\cdot$
  &
    $\hspace{-4mm} 10^{-14}\myi$
  &
    $\hspace{9mm} 3.6813\dotsc\cdot$
  &
    $\hspace{-4mm} 10^{-9}$
  \\ 
    $4$
  &
    $-3.7956\dotsc\cdot$
  &
    $\hspace{-4mm} 10^{-8}$
  &
    $\hspace{-4mm}+2.0698\dotsc\cdot$
  &
    $\hspace{-4mm} 10^{-7}\myi$
  &
    $\hspace{9mm} 1.7175\dotsc\cdot$
  &
    $\hspace{-4mm} 10^{-7}$
  \\ 
    $5$
  &
    $1.5216\dotsc\cdot$
  &
    $\hspace{-4mm} 10^{-3}$
  &
    $\hspace{-4mm}+8.8240\dotsc\cdot$
  &
    $\hspace{-4mm} 10^{-3}\myi$
  &
    $\hspace{9mm} 8.2908\dotsc\cdot$
  &
    $\hspace{-4mm} 10^{-6}$
  \\ 
    $6$
  &
    $3.7984\dotsc\cdot$
  &
    $\hspace{-4mm} 10^{0}$
  &
    $\hspace{-4mm}-4.2877\dotsc\cdot$
  &
    $\hspace{-4mm} 10^{0}\myi$
  &
    $\hspace{9mm} 1.0541\dotsc\cdot$
  &
    $\hspace{-4mm} 10^{-4}$
  \\ 
    $7$
  &
    $-3.0090\dotsc\cdot$
  &
    $\hspace{-4mm} 10^{2}$
  &
    $\hspace{-4mm}+4.8652\dotsc\cdot$
  &
    $\hspace{-4mm} 10^{2}\myi$
  &
    $\hspace{9mm} 2.3986\dotsc\cdot$
  &
    $\hspace{-4mm} 10^{-3}$
  \\ 
    $8$
  &
    $-1.8858\dotsc\cdot$
  &
    $\hspace{-4mm} 10^{3}$
  &
    $\hspace{-4mm}+4.2625\dotsc\cdot$
  &
    $\hspace{-4mm} 10^{2}\myi$
  &
    $\hspace{9mm} 1.8282\dotsc\cdot$
  &
    $\hspace{-4mm} 10^{-2}$
  \\ 
    $9$
  &
    $-1.1127\dotsc\cdot$
  &
    $\hspace{-4mm} 10^{3}$
  &
    $\hspace{-4mm}+9.0420\dotsc\cdot$
  &
    $\hspace{-4mm} 10^{2}\myi$
  &
    $\hspace{9mm} 5.7323\dotsc\cdot$
  &
    $\hspace{-4mm} 10^{-2}$
  \\ 
    $10$
  &
    $-5.2122\dotsc\cdot$
  &
    $\hspace{-4mm} 10^{1}$
  &
    $\hspace{-4mm}+2.1236\dotsc\cdot$
  &
    $\hspace{-4mm} 10^{1}\myi$
  &
    $\hspace{9mm} 1.0387\dotsc\cdot$
  &
    $\hspace{-4mm} 10^{0}$
  \\ 
    $11$
  &
    $1.0390\dotsc\cdot$
  &
    $\hspace{-4mm} 10^{1}$
  &
    $\hspace{-4mm}-7.6856\dotsc\cdot$
  &
    $\hspace{-4mm} 10^{0}\myi$
  &
    $\hspace{9mm} 5.2691\dotsc\cdot$
  &
    $\hspace{-4mm} 10^{-1}$
  \\ 
  \hhline{|b:=:b:======:b|}
\end{tabular}

%% file: circle_drop_new.tex
\begin{tabular}{||l||l||}
  \hhline{|t:=:t:=:t|}
    \multicolumn{1}{||c||}{$\M_{\D}$}&\multicolumn{1}{c||}{\vpad{$\vert Q(\mathfrak{A},\mathfrak{D})-1\vert$}}
  \\ 
  \hhline{|:=::=:|}
    $\{\}$
  &
    $1.4487\dotsc\cdot 10^{-26}$
  \\ 
    $\{1\}$
  &
    $4.5365\dotsc\cdot 10^{-25}$
  \\ 
    $\{2\}$
  &
    $2.1994\dotsc\cdot 10^{-25}$
  \\ 
    $\{5\}$
  &
    $8.0930\dotsc\cdot 10^{-26}$
  \\ 
    $\{10\}$
  &
    $3.5842\dotsc\cdot 10^{-26}$
  \\ 
    $\{15\}$
  &
    $2.1416\dotsc\cdot 10^{-26}$
  \\ 
    $\{1, 2\}$
  &
    $7.5493\dotsc\cdot 10^{-24}$
  \\ 
    $\{2, 3\}$
  &
    $2.3575\dotsc\cdot 10^{-24}$
  \\ 
    $\{2, 5\}$
  &
    $1.3344\dotsc\cdot 10^{-24}$
  \\ 
    $\{3, 7\}$
  &
    $5.8436\dotsc\cdot 10^{-25}$
  \\ 
    $\{2, 4, 6\}$
  &
    $9.1420\dotsc\cdot 10^{-24}$
  \\ 
    $\{5, 11, 14\}$
  &
    $3.5669\dotsc\cdot 10^{-25}$
  \\ 
    $\{1, 3, 5, 7\}$
  &
    $1.5529\dotsc\cdot 10^{-22}$
  \\ 
    $\{2, 4, 6, 8\}$
  &
    $3.7291\dotsc\cdot 10^{-23}$
  \\ 
    $\{3, 4, 5, 6, 7\}$
  &
    $2.1682\dotsc\cdot 10^{-22}$
  \\ 
    $\{3, 6, 9, 12, 15\}$
  &
    $1.0883\dotsc\cdot 10^{-23}$
  \\ 
    $\{9, 10, 11, 12, 13\}$
  &
    $1.5617\dotsc\cdot 10^{-24}$
  \\ 
    $\{2, 4, 6, 8, 10, 12\}$
  &
    $3.4899\dotsc\cdot 10^{-22}$
  \\ 
    $\{1, 3, 5, 7, 9, 11\}$
  &
    $1.8820\dotsc\cdot 10^{-21}$
  \\ 
    $\{2, 4, 6, 8, 10, 12, 14\}$
  &
    $8.5487\dotsc\cdot 10^{-22}$
  \\ 
    $\{2, 3, 4, 6, 8, 10, 12\}$
  &
    $5.6029\dotsc\cdot 10^{-21}$
  \\ 
  \hhline{|b:=:b:=:b|}
\end{tabular}

%% file: R3_new.tex
\begin{tabular}{||l||l|l|l||}
  \hhline{|t:=:t:===:t|}
    \multicolumn{1}{||c||}{}&\multicolumn{3}{c||}{\vpad{$\vert S_N(a)-1\vert$}}\\\hhline{||~||-|-|-||}\multicolumn{1}{||c||}{$a$}&\multicolumn{1}{c|}{$N=20$}&\multicolumn{1}{c|}{$N=70$}&\multicolumn{1}{c||}{$N=200$}
  \\ 
  \hhline{|:=::===:|}
    $-6$
  &
    $3.9870\dotsc\cdot 10^{-18}$
  &
    $1.7014\dotsc\cdot 10^{-105}$
  &
    $3.7560\dotsc\cdot 10^{-368}$
  \\ 
    $-3$
  &
    $4.7492\dotsc\cdot 10^{-22}$
  &
    $6.5289\dotsc\cdot 10^{-111}$
  &
    $6.8468\dotsc\cdot 10^{-375}$
  \\ 
    $-1$
  &
    $1.2054\dotsc\cdot 10^{-24}$
  &
    $1.6130\dotsc\cdot 10^{-114}$
  &
    $2.2057\dotsc\cdot 10^{-379}$
  \\ 
    $-1+300\mathrm{i}$
  &
    $2.8175\dotsc\cdot 10^{-16}$
  &
    $1.5955\dotsc\cdot 10^{-98}$
  &
    $1.0048\dotsc\cdot 10^{-378}$
  \\ 
    $-0.8$
  &
    $6.6429\dotsc\cdot 10^{-25}$
  &
    $7.0318\dotsc\cdot 10^{-115}$
  &
    $7.8416\dotsc\cdot 10^{-380}$
  \\ 
    $-0.8+50\mathrm{i}$
  &
    $4.9124\dotsc\cdot 10^{-24}$
  &
    $1.6522\dotsc\cdot 10^{-115}$
  &
    $4.8600\dotsc\cdot 10^{-380}$
  \\ 
    $0$
  &
    $6.1439\dotsc\cdot 10^{-26}$
  &
    $2.5409\dotsc\cdot 10^{-116}$
  &
    $1.2527\dotsc\cdot 10^{-381}$
  \\ 
    $50\mathrm{i}$
  &
    $5.1433\dotsc\cdot 10^{-25}$
  &
    $6.3145\dotsc\cdot 10^{-117}$
  &
    $7.8169\dotsc\cdot 10^{-382}$
  \\ 
    $300\mathrm{i}$
  &
    $1.4208\dotsc\cdot 10^{-17}$
  &
    $2.4131\dotsc\cdot 10^{-100}$
  &
    $6.6292\dotsc\cdot 10^{-381}$
  \\ 
    $1500\mathrm{i}$
  &
    $3.5498\dotsc\cdot 10^{-14}$
  &
    $3.8431\dotsc\cdot 10^{-80}$
  &
    $3.1958\dotsc\cdot 10^{-307}$
  \\ 
    $0.2+10\mathrm{i}$
  &
    $2.4939\dotsc\cdot 10^{-26}$
  &
    $1.0349\dotsc\cdot 10^{-116}$
  &
    $4.3685\dotsc\cdot 10^{-382}$
  \\ 
    $0.2+100\mathrm{i}$
  &
    $1.6846\dotsc\cdot 10^{-21}$
  &
    $3.3113\dotsc\cdot 10^{-117}$
  &
    $8.0859\dotsc\cdot 10^{-383}$
  \\ 
    $0.5$
  &
    $1.3906\dotsc\cdot 10^{-26}$
  &
    $3.1908\dotsc\cdot 10^{-117}$
  &
    $9.4420\dotsc\cdot 10^{-383}$
  \\ 
    $0.5+50\mathrm{i}$
  &
    $1.2520\dotsc\cdot 10^{-25}$
  &
    $8.2073\dotsc\cdot 10^{-118}$
  &
    $5.9166\dotsc\cdot 10^{-383}$
  \\ 
    $0.5+1500\mathrm{i}$
  &
    $6.7418\dotsc\cdot 10^{-15}$
  &
    $1.0285\dotsc\cdot 10^{-80}$
  &
    $4.0432\dotsc\cdot 10^{-308}$
  \\ 
    $0.8$
  &
    $5.7069\dotsc\cdot 10^{-27}$
  &
    $9.1900\dotsc\cdot 10^{-118}$
  &
    $2.0017\dotsc\cdot 10^{-383}$
  \\ 
    $0.8+300\mathrm{i}$
  &
    $1.1103\dotsc\cdot 10^{-18}$
  &
    $7.2488\dotsc\cdot 10^{-102}$
  &
    $1.1928\dotsc\cdot 10^{-382}$
  \\ 
    $0.8+1500\mathrm{i}$
  &
    $2.7648\dotsc\cdot 10^{-15}$
  &
    $4.2959\dotsc\cdot 10^{-81}$
  &
    $1.0755\dotsc\cdot 10^{-308}$
  \\ 
    $1+10\mathrm{i}$
  &
    $2.3862\dotsc\cdot 10^{-27}$
  &
    $3.7524\dotsc\cdot 10^{-118}$
  &
    $6.9829\dotsc\cdot 10^{-384}$
  \\ 
    $1+50\mathrm{i}$
  &
    $3.0417\dotsc\cdot 10^{-26}$
  &
    $1.0666\dotsc\cdot 10^{-118}$
  &
    $4.4787\dotsc\cdot 10^{-384}$
  \\ 
    $1+300\mathrm{i}$
  &
    $5.7549\dotsc\cdot 10^{-19}$
  &
    $2.9804\dotsc\cdot 10^{-102}$
  &
    $4.3682\dotsc\cdot 10^{-383}$
  \\ 
    $2+10\mathrm{i}$
  &
    $1.2736\dotsc\cdot 10^{-28}$
  &
    $5.9438\dotsc\cdot 10^{-120}$
  &
    $3.9701\dotsc\cdot 10^{-386}$
  \\ 
    $2+1500\mathrm{i}$
  &
    $1.3975\dotsc\cdot 10^{-16}$
  &
    $7.2077\dotsc\cdot 10^{-83}$
  &
    $1.7403\dotsc\cdot 10^{-311}$
  \\ 
  \hhline{|b:=:b:===:b|}
\end{tabular}

%% file: conj3new.tex
\begin{tabular}{||l||l|l|l||}
  \hhline{|t:=:t:===:t|}
    \multicolumn{1}{||c||}{}&\multicolumn{3}{c||}{\vpad{$\max_{1\le n\le m+1}{\left\vert\frac{n-m-2}{n\,\eta^{\langle l+m \rangle}(a)}\cdot\frac{E_{l,m,N,n-1}(a)}{E_{l,m,N,n}(a)}
 -1\right\vert}$}}\\\hhline{||~||-|-|-||}\multicolumn{1}{||c||}{$m$}&\multicolumn{1}{c|}{$l=2$}&\multicolumn{1}{c|}{$l=5$}&\multicolumn{1}{c||}{$l=8$}
  \\ 
  \hhline{|:=::===:|}
    $1$
  &
    $8.336\dotsc\cdot 10^{-75}$
  &
    $9.685\dotsc\cdot 10^{-74}$
  &
    $1.023\dotsc\cdot 10^{-72}$
  \\ 
    $2$
  &
    $6.211\dotsc\cdot 10^{-71}$
  &
    $6.871\dotsc\cdot 10^{-70}$
  &
    $6.949\dotsc\cdot 10^{-69}$
  \\ 
    $3$
  &
    $3.266\dotsc\cdot 10^{-67}$
  &
    $3.455\dotsc\cdot 10^{-66}$
  &
    $3.342\dotsc\cdot 10^{-65}$
  \\ 
    $4$
  &
    $1.296\dotsc\cdot 10^{-63}$
  &
    $1.312\dotsc\cdot 10^{-62}$
  &
    $1.213\dotsc\cdot 10^{-61}$
  \\ 
    $5$
  &
    $4.044\dotsc\cdot 10^{-60}$
  &
    $3.914\dotsc\cdot 10^{-59}$
  &
    $3.461\dotsc\cdot 10^{-58}$
  \\ 
    $6$
  &
    $1.019\dotsc\cdot 10^{-56}$
  &
    $9.426\dotsc\cdot 10^{-56}$
  &
    $7.968\dotsc\cdot 10^{-55}$
  \\ 
    $7$
  &
    $2.114\dotsc\cdot 10^{-53}$
  &
    $1.868\dotsc\cdot 10^{-52}$
  &
    $1.510\dotsc\cdot 10^{-51}$
  \\ 
    $8$
  &
    $3.664\dotsc\cdot 10^{-50}$
  &
    $3.093\dotsc\cdot 10^{-49}$
  &
    $2.391\dotsc\cdot 10^{-48}$
  \\ 
  \hhline{|b:=:b:===:b|}
\end{tabular}

%% file: R3fact_new.tex
\begin{tabular}{||l||l|l|l||}
  \hhline{|t:=:t:===:t|}
    \multicolumn{1}{||c||}{}&\multicolumn{3}{c||}{\vpad{$\vert S_{!\,N}(a)-1\vert$}}\\\hhline{||~||-|-|-||}\multicolumn{1}{||c||}{$a$}&\multicolumn{1}{c|}{$N=70$}&\multicolumn{1}{c|}{$N=200$}&\multicolumn{1}{c||}{$N=700$}
  \\ 
  \hhline{|:=::===:|}
    $-2$
  &
    $6.6474\dotsc\cdot 10^{-3}$
  &
    $4.2041\dotsc\cdot 10^{-5}$
  &
    $1.0842\dotsc\cdot 10^{-13}$
  \\ 
    $-1$
  &
    $9.9362\dotsc\cdot 10^{-3}$
  &
    $8.6712\dotsc\cdot 10^{-6}$
  &
    $1.8735\dotsc\cdot 10^{-14}$
  \\ 
    $-0.8$
  &
    $1.1821\dotsc\cdot 10^{-2}$
  &
    $1.0637\dotsc\cdot 10^{-5}$
  &
    $1.2612\dotsc\cdot 10^{-14}$
  \\ 
    $0$
  &
    $1.2694\dotsc\cdot 10^{-2}$
  &
    $1.5834\dotsc\cdot 10^{-6}$
  &
    $2.3561\dotsc\cdot 10^{-15}$
  \\ 
    $50\myi$
  &
    $9.9046\dotsc\cdot 10^{-2}$
  &
    $5.6900\dotsc\cdot 10^{-5}$
  &
    $1.0808\dotsc\cdot 10^{-14}$
  \\ 
    $0.2$
  &
    $3.9934\dotsc\cdot 10^{-3}$
  &
    $8.0428\dotsc\cdot 10^{-7}$
  &
    $3.7177\dotsc\cdot 10^{-14}$
  \\ 
    $0.2+10\myi$
  &
    $9.6420\dotsc\cdot 10^{-4}$
  &
    $1.1453\dotsc\cdot 10^{-6}$
  &
    $4.8160\dotsc\cdot 10^{-15}$
  \\ 
    $0.5$
  &
    $1.2866\dotsc\cdot 10^{-3}$
  &
    $3.1413\dotsc\cdot 10^{-7}$
  &
    $8.2725\dotsc\cdot 10^{-16}$
  \\ 
    $0.8$
  &
    $4.8940\dotsc\cdot 10^{-4}$
  &
    $1.6870\dotsc\cdot 10^{-7}$
  &
    $3.1938\dotsc\cdot 10^{-16}$
  \\ 
    $0.8+10\myi$
  &
    $2.5509\dotsc\cdot 10^{-4}$
  &
    $2.6686\dotsc\cdot 10^{-7}$
  &
    $7.1105\dotsc\cdot 10^{-16}$
  \\ 
    $1+10\myi$
  &
    $3.4554\dotsc\cdot 10^{-4}$
  &
    $2.3766\dotsc\cdot 10^{-7}$
  &
    $4.4841\dotsc\cdot 10^{-16}$
  \\ 
    $1+50\myi$
  &
    $8.7300\dotsc\cdot 10^{-2}$
  &
    $4.6678\dotsc\cdot 10^{-6}$
  &
    $4.5817\dotsc\cdot 10^{-16}$
  \\ 
    $2+10\myi$
  &
    $2.4417\dotsc\cdot 10^{-4}$
  &
    $1.4395\dotsc\cdot 10^{-7}$
  &
    $2.5638\dotsc\cdot 10^{-17}$
  \\ 
  \hhline{|b:=:b:===:b|}
\end{tabular}

%% file: Rfact3.tex
\begin{tabular}{||l||l|l|l||}
  \hhline{|t:=:t:===:t|}
    \multicolumn{1}{||c||}{}&\multicolumn{3}{c||}{\vpad{$\vert S^{!}_{N}(a)-1\vert$}}\\\hhline{||~||-|-|-||}\multicolumn{1}{||c||}{$a$}&\multicolumn{1}{c|}{$N=70$}&\multicolumn{1}{c|}{$N=200$}&\multicolumn{1}{c||}{$N=700$}
  \\ 
  \hhline{|:=::===:|}
    $-2$
  &
    $2.2873\dotsc\cdot 10^{-1}$
  &
    $2.9563\dotsc\cdot 10^{-2}$
  &
    $1.0597\dotsc\cdot 10^{-2}$
  \\ 
    $-1$
  &
    $4.0733\dotsc\cdot 10^{-2}$
  &
    $1.8496\dotsc\cdot 10^{-2}$
  &
    $6.7039\dotsc\cdot 10^{-3}$
  \\ 
    $-0.8$
  &
    $3.6069\dotsc\cdot 10^{-2}$
  &
    $2.0192\dotsc\cdot 10^{-2}$
  &
    $5.9950\dotsc\cdot 10^{-3}$
  \\ 
    $0$
  &
    $2.3313\dotsc\cdot 10^{-2}$
  &
    $1.1531\dotsc\cdot 10^{-2}$
  &
    $4.1362\dotsc\cdot 10^{-3}$
  \\ 
    $50\myi$
  &
    $7.2888\dotsc\cdot 10^{-2}$
  &
    $2.7648\dotsc\cdot 10^{-2}$
  &
    $9.1542\dotsc\cdot 10^{-3}$
  \\ 
    $0.2$
  &
    $2.0458\dotsc\cdot 10^{-2}$
  &
    $6.6277\dotsc\cdot 10^{-3}$
  &
    $3.7877\dotsc\cdot 10^{-3}$
  \\ 
    $0.2+10\myi$
  &
    $6.8470\dotsc\cdot 10^{-2}$
  &
    $3.0200\dotsc\cdot 10^{-2}$
  &
    $1.1262\dotsc\cdot 10^{-2}$
  \\ 
    $0.5$
  &
    $1.5076\dotsc\cdot 10^{-2}$
  &
    $7.6916\dotsc\cdot 10^{-3}$
  &
    $3.0248\dotsc\cdot 10^{-3}$
  \\ 
    $0.8$
  &
    $1.3358\dotsc\cdot 10^{-2}$
  &
    $6.7687\dotsc\cdot 10^{-3}$
  &
    $6.9181\dotsc\cdot 10^{-4}$
  \\ 
    $0.8+10\myi$
  &
    $3.8167\dotsc\cdot 10^{-2}$
  &
    $1.7013\dotsc\cdot 10^{-2}$
  &
    $6.3566\dotsc\cdot 10^{-3}$
  \\ 
    $1+50\myi$
  &
    $2.1107\dotsc\cdot 10^{-2}$
  &
    $8.5441\dotsc\cdot 10^{-3}$
  &
    $2.9618\dotsc\cdot 10^{-3}$
  \\ 
    $2+10\myi$
  &
    $1.2910\dotsc\cdot 10^{-2}$
  &
    $5.7831\dotsc\cdot 10^{-3}$
  &
    $2.1859\dotsc\cdot 10^{-3}$
  \\ 
  \hhline{|b:=:b:===:b|}
\end{tabular}

%% file: 1stmissing.tex
\begin{tabular}{||l|l||l|l||}
    \hhline{|t:==:t:==:t|}
    \multicolumn{1}{||c|}{}&\multicolumn{1}{c||}{}&\multicolumn{2}{c||}{$S_{N,\mathfrak{D},\mathfrak{F}}(a)$}
  \\ 
    \hhline{||~|~||-|-||}\multicolumn{1}{||c|}{$\mathfrak{M}_\mathfrak{F}$}&\multicolumn{1}{c||}{$\mathfrak{M}_\mathfrak{D}$}&\multicolumn{1}{c|}{$N=40$}&\multicolumn{1}{c||}{$N=150$}
  \\ 
    \hhline{|:==::==:|}
    $\{1\}$
  &
    $\{1\}$
  &
    $1.5705\dotsc\cdot 10^{-57}$
  &
    $4.1503\dotsc\cdot 10^{-273}$
  \\ 
    $$
  &
    $\{2, 3\}$
  &
    $1.5883\dotsc\cdot 10^{-56}$
  &
    $1.4302\dotsc\cdot 10^{-271}$
  \\ 
    $$
  &
    $\{1, 2, 3\}$
  &
    $1.0588\dotsc\cdot 10^{-54}$
  &
    $3.0229\dotsc\cdot 10^{-269}$
  \\ 
    $$
  &
    $\{4, 5\}$
  &
    $4.4938\dotsc\cdot 10^{-57}$
  &
    $4.2342\dotsc\cdot 10^{-272}$
  \\ 
    $$
  &
    $\{2, 3, 5\}$
  &
    $1.9855\dotsc\cdot 10^{-55}$
  &
    $5.9628\dotsc\cdot 10^{-270}$
  \\ 
    $$
  &
    $\{1, 3, 4, 5, 7\}$
  &
    $6.3476\dotsc\cdot 10^{-53}$
  &
    $1.9159\dotsc\cdot 10^{-266}$
  \\ 
    \hhline{||-|-||-|-||}
    $\{2, 3\}$
  &
    $\{1\}$
  &
    $1.5883\dotsc\cdot 10^{-56}$
  &
    $1.4302\dotsc\cdot 10^{-271}$
  \\ 
    $$
  &
    $\{2, 3\}$
  &
    $1.6076\dotsc\cdot 10^{-55}$
  &
    $4.9296\dotsc\cdot 10^{-270}$
  \\ 
    $$
  &
    $\{1, 2, 3\}$
  &
    $1.0726\dotsc\cdot 10^{-53}$
  &
    $1.0420\dotsc\cdot 10^{-267}$
  \\ 
    $$
  &
    $\{4, 5\}$
  &
    $4.5482\dotsc\cdot 10^{-56}$
  &
    $1.4593\dotsc\cdot 10^{-270}$
  \\ 
    $$
  &
    $\{2, 3, 5\}$
  &
    $2.0112\dotsc\cdot 10^{-54}$
  &
    $2.0554\dotsc\cdot 10^{-268}$
  \\ 
    $$
  &
    $\{1, 3, 4, 5, 7\}$
  &
    $6.4409\dotsc\cdot 10^{-52}$
  &
    $6.6060\dotsc\cdot 10^{-265}$
  \\ 
    \hhline{||-|-||-|-||}
    $\{1, 2, 3\}$
  &
    $\{1\}$
  &
    $1.0588\dotsc\cdot 10^{-54}$
  &
    $3.0229\dotsc\cdot 10^{-269}$
  \\ 
    $$
  &
    $\{2, 3\}$
  &
    $1.0726\dotsc\cdot 10^{-53}$
  &
    $1.0420\dotsc\cdot 10^{-267}$
  \\ 
    $$
  &
    $\{1, 2, 3\}$
  &
    $7.1626\dotsc\cdot 10^{-52}$
  &
    $2.2028\dotsc\cdot 10^{-265}$
  \\ 
    $$
  &
    $\{4, 5\}$
  &
    $3.0343\dotsc\cdot 10^{-54}$
  &
    $3.0848\dotsc\cdot 10^{-268}$
  \\ 
    $$
  &
    $\{2, 3, 5\}$
  &
    $1.3430\dotsc\cdot 10^{-52}$
  &
    $4.3452\dotsc\cdot 10^{-266}$
  \\ 
    $$
  &
    $\{1, 3, 4, 5, 7\}$
  &
    $4.3085\dotsc\cdot 10^{-50}$
  &
    $1.3968\dotsc\cdot 10^{-262}$
  \\ 
    \hhline{||-|-||-|-||}
    $\{4, 5\}$
  &
    $\{1\}$
  &
    $4.4938\dotsc\cdot 10^{-57}$
  &
    $4.2342\dotsc\cdot 10^{-272}$
  \\ 
    $$
  &
    $\{2, 3\}$
  &
    $4.5482\dotsc\cdot 10^{-56}$
  &
    $1.4593\dotsc\cdot 10^{-270}$
  \\ 
    $$
  &
    $\{1, 2, 3\}$
  &
    $3.0343\dotsc\cdot 10^{-54}$
  &
    $3.0848\dotsc\cdot 10^{-268}$
  \\ 
    $$
  &
    $\{4, 5\}$
  &
    $1.2867\dotsc\cdot 10^{-56}$
  &
    $4.3204\dotsc\cdot 10^{-271}$
  \\ 
    $$
  &
    $\{2, 3, 5\}$
  &
    $5.6898\dotsc\cdot 10^{-55}$
  &
    $6.0849\dotsc\cdot 10^{-269}$
  \\ 
    $$
  &
    $\{1, 3, 4, 5, 7\}$
  &
    $1.8219\dotsc\cdot 10^{-52}$
  &
    $1.9556\dotsc\cdot 10^{-265}$
  \\ 
    \hhline{||-|-||-|-||}
    $\{2, 3, 5\}$
  &
    $\{1\}$
  &
    $1.9855\dotsc\cdot 10^{-55}$
  &
    $5.9628\dotsc\cdot 10^{-270}$
  \\ 
    $$
  &
    $\{2, 3\}$
  &
    $2.0112\dotsc\cdot 10^{-54}$
  &
    $2.0554\dotsc\cdot 10^{-268}$
  \\ 
    $$
  &
    $\{1, 2, 3\}$
  &
    $1.3430\dotsc\cdot 10^{-52}$
  &
    $4.3452\dotsc\cdot 10^{-266}$
  \\ 
    $$
  &
    $\{4, 5\}$
  &
    $5.6898\dotsc\cdot 10^{-55}$
  &
    $6.0849\dotsc\cdot 10^{-269}$
  \\ 
    $$
  &
    $\{2, 3, 5\}$
  &
    $2.5183\dotsc\cdot 10^{-53}$
  &
    $8.5712\dotsc\cdot 10^{-267}$
  \\ 
    $$
  &
    $\{1, 3, 4, 5, 7\}$
  &
    $8.0783\dotsc\cdot 10^{-51}$
  &
    $2.7554\dotsc\cdot 10^{-263}$
  \\ 
    \hhline{||-|-||-|-||}
    $\{1, 3, 4, 5, 7\}$
  &
    $\{1\}$
  &
    $6.3476\dotsc\cdot 10^{-53}$
  &
    $1.9159\dotsc\cdot 10^{-266}$
  \\ 
    $$
  &
    $\{2, 3\}$
  &
    $6.4409\dotsc\cdot 10^{-52}$
  &
    $6.6060\dotsc\cdot 10^{-265}$
  \\ 
    $$
  &
    $\{1, 2, 3\}$
  &
    $4.3085\dotsc\cdot 10^{-50}$
  &
    $1.3968\dotsc\cdot 10^{-262}$
  \\ 
    $$
  &
    $\{4, 5\}$
  &
    $1.8219\dotsc\cdot 10^{-52}$
  &
    $1.9556\dotsc\cdot 10^{-265}$
  \\ 
    $$
  &
    $\{2, 3, 5\}$
  &
    $8.0783\dotsc\cdot 10^{-51}$
  &
    $2.7554\dotsc\cdot 10^{-263}$
  \\ 
    $$
  &
    $\{1, 3, 4, 5, 7\}$
  &
    $2.6003\dotsc\cdot 10^{-48}$
  &
    $8.8624\dotsc\cdot 10^{-260}$
  \\ 
    \hhline{|b:==:b:==:b|}
\end{tabular}

%% file: gen_char_poly.tex
\begin{tabular}{||l||l|l|l||}
  \hhline{|t:=:t:===:t|}
    \multicolumn{1}{||c||}{}&\multicolumn{3}{c||}{\vpad{$\left\vert \frac{V_{l,n,\alpha}(a)}{\eta^{\langle l \rangle}(a)W_{l,n,\alpha}(a)}-1\right\vert$}}\\\hhline{||~||-|-|-||}\multicolumn{1}{||c||}{$n$}&\multicolumn{1}{c|}{$\alpha=\alpha_{17,3}$}&\multicolumn{1}{c|}{$\alpha=\alpha_{19,3}$}&\multicolumn{1}{c||}{$\alpha=\alpha_{23,3}$}
  \\ 
  \hhline{|:=::===:|}
    $0$
  &
    $4.981\dotsc\cdot 10^{-22}$
  &
    $3.585\dotsc\cdot 10^{-25}$
  &
    $1.236\dotsc\cdot 10^{-31}$
  \\ 
    $1$
  &
    $2.100\dotsc\cdot 10^{-19}$
  &
    $1.721\dotsc\cdot 10^{-22}$
  &
    $7.387\dotsc\cdot 10^{-29}$
  \\ 
    $2$
  &
    $4.741\dotsc\cdot 10^{-17}$
  &
    $4.302\dotsc\cdot 10^{-20}$
  &
    $2.314\dotsc\cdot 10^{-26}$
  \\ 
    $3$
  &
    $2.157\dotsc\cdot 10^{-14}$
  &
    $3.209\dotsc\cdot 10^{-17}$
  &
    $2.695\dotsc\cdot 10^{-23}$
  \\ 
    $4$
  &
    $6.251\dotsc\cdot 10^{-12}$
  &
    $2.914\dotsc\cdot 10^{-14}$
  &
    $6.718\dotsc\cdot 10^{-21}$
  \\ 
    $5$
  &
    $1.049\dotsc\cdot 10^{-9}$
  &
    $2.096\dotsc\cdot 10^{-12}$
  &
    $8.715\dotsc\cdot 10^{-19}$
  \\ 
    $6$
  &
    $7.331\dotsc\cdot 10^{-8}$
  &
    $5.602\dotsc\cdot 10^{-11}$
  &
    $9.066\dotsc\cdot 10^{-16}$
  \\ 
    $7$
  &
    $6.999\dotsc\cdot 10^{-7}$
  &
    $1.894\dotsc\cdot 10^{-8}$
  &
    $2.404\dotsc\cdot 10^{-13}$
  \\ 
    $8$
  &
    $6.011\dotsc\cdot 10^{-4}$
  &
    $1.650\dotsc\cdot 10^{-6}$
  &
    $5.683\dotsc\cdot 10^{-11}$
  \\ 
  \hhline{|b:=:b:===:b|}
\end{tabular}